\documentclass[12pt,a4paper]{article}

\usepackage{amssymb}

\usepackage{amsmath}

\newtheorem{theorem}{Theorem}[section]
\newtheorem{definition}[theorem]{Definition}
\newtheorem{lemma}[theorem]{Lemma}
\newtheorem{corollary}[theorem]{Corollary}
\newtheorem{proposition}[theorem]{Proposition}

\newtheorem{algorithm}[theorem]{Algorithm}

\begin{document}
\title{Gr\"{o}bner-Shirshov bases for metabelian Lie algebras\footnote{Supported by the NNSF of China (Nos. 10771077;
10911120389).}}
\author{ Yongshan Chen and
 Yuqun
Chen\footnote{Corresponding author.}
\\
{\small \ School of Mathematical Sciences, South China Normal
University}\\
{\small Guangzhou 510631, P. R. China}\\
{\small jackalshan@126.com}\\
{\small Email: yqchen@scnu.edu.cn}}

\date{}

\maketitle \noindent\textbf{Abstract:}  In this paper, we establish
the Gr\"{o}bner-Shirshov bases theory for metabelian Lie algebras.
As applications, we find the Gr\"{o}bner-Shirshov bases for partial
commutative metabelian Lie algebras related to circuits, trees and
some cubes.



\section{Introduction}
The class of metabelian Lie algebras is an important class of Lie
algebras and attracts many attentions. Let us mention the recent
papers by E. Daniyarova, I. Kazatchkov, and V. Remeslennikov
\cite{dkr1,dkr2,dkr3} on algebraic geometry of free metabelian Lie
algebra, S. Findik and V. Drensky \cite{f,df} on automorphisms of
free metabelian Lie algebras, and V. Kurlin \cite{k} on the
Backer-Campbell-Hausdorff formula for free metabelian Lie algebras.
Gr\"{o}bner-Shirshov bases theory would be useful on this class of
algebras. This theory was first considered by V.V. Talapov \cite{ta}
in 1982. However, there are serious gaps in his paper. He missed
several cases when he defined compositions. This means the theory
was not established correctly. We refine his idea and complete the
results.

It is well-known that for many kinds of algebras, if
$A_i=(X_i|S_i)$, $i=1,2$, are defined by generators and defining
relations, where $S_1$ and $S_2$ are Gr\"{o}bner-Shirshov bases
respectively, then $S_1\cup S_2$ is a Gr\"{o}bner-Shirshov basis for
the free product $A_1\ast A_2=(X_1\cup X_2|S_1\cup S_2)$ of $A_1$
and $A_2$, for example, associative algebras, Lie algebras and for
all classes with compositions of inclusion and intersection only
(cf. \cite{bc,bcs}). We prove that it is not the case for metabelian
Lie algebras, see Theorem \ref{t4.1}, even in the case of
$S_2=\varnothing$. On the other hand, if $S_i\subset A_i^{(2)}$,
then $S_1\cup S_2$ is a Gr\"{o}bner-Shirshov basis for the free
metabelian Lie product $A_1\ast A_2$, see Proposition \ref{p4.2}.

Throughout this paper, all algebras will be considered over a field
$\bf{k}$ of arbitrary characteristic. Suppose that $\mathcal{L}$ is
a Lie algebra. Then $\mathcal{L}$ is called a metabelian Lie algebra
if $\mathcal{L}^{(2)}=0$, where $\mathcal{L}^{(0)}=\mathcal{L}$,
$\mathcal{L}^{(n+1)}=[\mathcal{L}^{(n)}, \mathcal{L}^{(n)}]$. More
precisely, the variety of metabelian Lie algebras is given by the
identity
 $$
 (x_1x_2)(x_3x_4)=0.
 $$

\section{Composition-Diamond lemma for  metabelian Lie algebras}
Let us begin with the construction of a free metabelian Lie algebra.
Let $X$ be a set and $Lie(X)$ be the free Lie algebra generated by
$X$. Then $\mathcal{L}_{(2)}(X)=Lie(X)/Lie(X)^{(2)}$ is the free
metabelian Lie algebra generated by $X$. Any metabelian Lie algebra
$\mathcal{ML}$ is a homomorphic image of a free metabelian Lie
algebra generated by some $X$, that is, $\mathcal{ML}$ can be
presented by generators $X$ and defining relations $S$:
$\mathcal{ML}=\mathcal{L}_{(2)}(X  | S)$.

We call a non-associative monomial on $X$ is \emph{left-normed} if
it is of the form $(\cdots((ab)c)\cdots)d$. In the sequel, the
brackets in the expression of left-normed monomials are omitted.

Let $X$ be well-ordered.  For an arbitrary set of indices
$j_1,j_2,\cdots,j_m$, define an associative word
$$
\langle a_{j_1}\cdots a_{j_m}\rangle=a_{i_1}\cdots a_{i_m},
$$
where $ a_{i_1}\leq \cdots \leq a_{i_m} $ and $i_1,i_2,\cdots,i_m$ is a
permutation of the indices $j_1,j_2,\cdots,j_m$.

Let
$$
R=\{u=a_0a_1a_2\cdots a_n \ |  a_i\in X \  (0\leq i\leq n),  \
a_0>a_1\leq \cdots \leq a_n ,
 n\geq1\}
$$
and $N=X\cup R$, where $u=a_0a_1a_2\cdots a_n$ is left-normed.

Then $N$ forms a linear basis of the free metabelian Lie algebra
$\mathcal{L}_{(2)}(X)$, i.e., $\mathcal{L}_{(2)}(X)={\bf k}N$, see
\cite{b63}.

We call elements of $N$ regular words on $X$ and those of $R$
regular $R$-words. Therefore, for any $f\in \mathcal{L}_{(2)}(X)$,
$f$ has a unique presentation $f=f^{(1)}+f^{(0)}$, where $f^{(1)}\in
{\bf k}R$ and $f^{(0)}\in {\bf k}X$. Moreover, the multiplication
table of regular words is the following, $u\cdot v=0$ if both
$u,v\in R$, and
$$
a_0a_1a_2\cdots a_n\cdot b=\left\{
\begin{array}{l}
a_0\langle a_1a_2\cdots a_nb\rangle \ \ \ \ \ \ \ \ \ \ \ \ \ \ \ \ \ \ \ \ \ \  \ \ \  \mbox{if} \ a_1\leq b ,\\
a_0ba_1a_2\cdots a_n-a_1b\langle a_0a_2\cdots a_n\rangle \ \ \ \ \mbox{if}  \ a_1> b .
\end{array}\right.
$$

If $u=a_0a_1\cdots a_n\in R$, then the regular words $a_i \ (0\leq i\leq
n)$, $a_0\langle a_{i_1}\cdots a_{i_l}\rangle$ ($l\leq n, \
a_{i_1},\cdots, a_{i_l}$ is a subsequence of the sequence $a_1,\cdots,
a_n$) are called subwords of $u$. The words $a_i \ (2\leq
i\leq n)$, and also $a_1$ if $a_0>a_2$ are called strict subwords of
$u$.

Define the length of regular words:
$$
|a|=1,  \ |a_0a_1a_2\cdots a_n|=n+1,
$$
where $a,\ a_0,\dots,a_n\in X$. Now we order the set $N$
degree-lexicographically, i.e., for any $u,v\in N$,
$$
u>v  \ \mbox{if} \ |u|>|v|  \ \mbox{or} \  |u|=|v|, \  u>_{lex}v.
$$

Through out this paper, we will use this ordering.

The largest monomial occurring in $f\in \mathcal{L}_{(2)}(X)$ with
nonzero coefficient is called the leading word of $f$ and is denoted
by $\bar{f}$. Then we have $\overline{a_0a_1a_2\cdots a_n\cdot
b}=a_0\langle a_1a_2\cdots a_nb\rangle$ and $|\overline{u\cdot
b}|=|u|+1$. For any $f\in \mathcal{L}_{(2)}(X)$, we called $f$ to be
monic, $(1)$-monic and $(0)$-monic if the coefficients of $\bar{f}$,
$\overline{f^{(1)}}$ and $\overline{f^{(0)}}$ are $1$ respectively.

\begin{lemma}\label{l3.1}For any $u,v\in N$, if $u>v$ then
 $$(\forall b\in N)\ u\cdot b\neq0\Rightarrow \overline{u\cdot b}>\overline{v\cdot b}.$$
\end{lemma}
{\bf Proof.} The result is obvious if either $u,v\in X$ or
$|u|>|v|$. Suppose that $u=a_0a_1a_2\cdots a_n$,
$v=a'_0a'_1a'_2\cdots a'_n\in R$ and $b\in X$.  If $a_0>a'_0$ then
we are done. If $a_0=a'_0$, then
 $\langle a_1a_2\cdots a_nb\rangle>\langle a'_1a'_2\cdots a'_nb\rangle$ in $[X]$ since the deg-lex
 ordering on $[X]$ is monomial, where $[X]$ is the free
 commutative momoid generated by $X$. Now, the result follows. $\Box$

\ \

Let $S\subset \mathcal{L}_{(2)}(X)$. We denote $u_s=sv_1v_2\cdots
v_n$, where $v_i\in N, \ s\in S $ and $n\geq0 $. We call $u_s$ an $s$-word (or
$S$-word). It is clear that each element of the ideal $Id(S)$ of
$\mathcal{L}_{(2)}(X)$ generated by $S$ is a linear combination of
$S$-words.

\begin{definition}\label{d3.1}
Let $S\subset \mathcal{L}_{(2)}(X)$. Then the following two kinds of
polynomials are called normal $S$-words:
\begin{enumerate}
\item[(i)] $sa_1a_2\cdots a_n$, where $a_i\in X\ (1\leq i\leq n),\ a_1\leq a_2\leq \cdots \leq a_n$,
$s\in S$, $\bar{s}\neq a_1$ and $n\geq 0$;
\item[(ii)] $us$, where $u\in R$, $s\in S$ and $\bar{s}\neq u$.
\end{enumerate}
\end{definition}

By a simple observation, we have
$$
\overline{sa_1a_2\cdots a_n}=\left\{\begin{array}{l}
c_0\langle c_1\cdots c_ka_1a_2\cdots a_n\rangle \ \ \ \ \ \ \ \ \ \ \mbox{if} \ \bar{s}=c_0c_1\cdots c_k,\\
c_0a_1a_2\cdots a_n \ \ \ \ \ \ \ \ \ \ \ \ \ \ \ \ \ \ \ \ \ \ \mbox{if} \ \bar{s}=c_0>a_1,\\
a_1c_0a_2\cdots a_n \ \ \ \ \ \ \ \ \ \ \ \ \ \ \ \ \ \ \ \ \ \ \mbox{if} \ \bar{s}=c_0<a_1,
\end{array}\right.
$$
and $ \overline{us}=a_0\langle a_1\cdots
a_k\overline{s^{(0)}}\rangle, $ where $u=a_0\langle a_1\cdots
a_k\rangle$. That is to say, if $u_s$ is a normal $s$-word, then
$\overline{u_s}$ either contains $\bar{s}$ as a subword or contains
$\overline{s^{(0)}}$ as a strict subword.

 A regular
 word $u$ is called $S$-irreducible if for any $s\in S$, $u$ contains
  neither $\bar{s}$ as a subword nor $\overline{s^{(0)}}$ as a strict subword.
   Denote $Irr(S)$ the set of all $S$-irreducible words. This means
$$
Irr(S)=\{u \ | \ u \in N, u\neq\overline{v_s} \ \mbox{for any normal S-word}  \ v_s\}.
$$

\noindent {\bf Remark:} For any $s\in \mathcal{L}_{(2)}(X)$,
$$
sa_1a_2\cdots a_n=sa_1a_{j_2}\cdots a_{j_n},
$$
where $\langle a_{j_2}\cdots a_{j_n}\rangle=a_2\cdots a_n$.

\begin{lemma}\label{l3.2} Let $S\subset \mathcal{L}_{(2)}(X)$ and $Id(S)$ be the ideal of $\mathcal{L}_{(2)}(X)$
generated by $S$. Then for any $f\in Id(S)$,
$f$ can be written as a linear combination of normal $S$-words.
\end{lemma}
{\bf Proof.} It is suffice to show that any $S$-word
$u_s=su_1u_2\cdots u_n$ is a linear combination of normal $S$-words,
where $u_i\in N, \ 1\leq i\leq n$. We may assume that $s$ is monic.
The proof will be proceeded by induction on $n$.

There is nothing to prove if $n=0$.

Assume that $n=1$. If $\bar{s}\neq u_1$, then either $su_1$ or
$u_1s$ is normal.
 If $\bar{s}=u_1$, then $s=u_1+\sum_{\bar{s}>v_j\in N}\alpha_jv_j$, $\alpha_j\in {\bf k}$ and
$$su_1=s(s-\sum_{v_j<\bar{s}}\alpha_i v_j)=-s\sum_{v_j<\bar{s}}\alpha_jv_j=-\sum_{v_j<\bar{s}}\alpha_js v_j,
$$
where for each $j$, either $v_js$ or $sv_j$ is normal.

For $n\geq2$, if $\exists u_i\in R\ (i\geq2)$, then $su_1u_2\cdots
u_n=0$;  if $u_1\in R$, then $(su_1)a_2\cdots a_n=s(u_1a_2\cdots
a_n)$ which is the above case. So we may assume that
$u_s=sa_1a_2\cdots a_n$ is normal and $u_{n+1}=a\in X$. Then
$$
u_s\cdot u_{n+1}= sa_1a_2\cdots a_n\cdot a=sa_1\langle a_2\cdots a_na\rangle.
$$
If $a\geq a_1$, then $sa_1\langle a_2\cdots a_na\rangle$ is normal. If $a<a_1$, then
\begin{eqnarray*}
u_s\cdot u_{n+1}&=&sa_1aa_2\cdots a_n\\
&=&saa_1a_2\cdots a_n-((a_1a)s)a_2\cdots a_n\\
&=&saa_1a_2\cdots a_n-a_1aa_2\cdots a_n\cdot s.
\end{eqnarray*}
Clearly, by the previous proof, $a_1aa_2\cdots a_n\cdot s$ is
normal. Now $saa_1a_2\cdots a_n$ is already normal provided that
$\bar{s}\neq a$. If $\bar{s}=a$, then we substitute $a$ by
$-\sum_{\bar{s}>v_j\in N}\alpha_jv_j$  where
$s=a+\sum_{\bar{s}>v_j\in N}\alpha_jv_j$, and the result follows
now. $\Box$

\begin{lemma}\label{l3} Let $u_s$ be a normal $S$-word and $w\in N$. If $\overline{u_s}<w$, then
$$
(\forall a\in X) \ w\cdot a\neq 0\Rightarrow \overline{u_s\cdot
a}<\overline{w\cdot a}.
$$
\end{lemma}
{\bf Proof.} Suppose that $w=b_0b_1\cdots b_m$ where $m\geq0$. Then
$$
\overline{w\cdot a}=\left\{\begin{array}{l}
b_0\langle b_1\cdots b_ma\rangle \ \ \ \ \ \ \mbox{if} \ \ m>0,\\
b_0a \ \ \ \ \ \ \ \ \ \ \ \ \ \  \ \  \ \ \ \mbox{if} \ \ m=0 \ \mbox{and} \ b_0>a,\\
ab_0 \ \ \ \ \ \ \ \ \ \ \ \ \ \  \ \ \ \ \ \mbox{if} \ \ m=0 \ \mbox{and} \ b_0<a.\\
\end{array}\right.
$$

If $u_s=sa_1a_2\cdots a_n$, then
$$
\overline{u_s}=\left\{\begin{array}{l}
c_0\langle c_1\cdots c_ka_1a_2\cdots a_n\rangle \ \ \ \ \ \ \ \ \ \ \mbox{if} \ \bar{s}=c_0c_1\cdots c_k,\\
c_0a_1a_2\cdots a_n \ \ \ \ \ \ \ \ \ \ \ \ \ \ \ \ \ \ \ \ \ \ \mbox{if} \ \bar{s}=c_0>a_1,\\
a_1c_0a_2\cdots a_n \ \ \ \ \ \ \ \ \ \ \ \ \ \ \ \ \ \ \ \ \ \
\mbox{if} \ \bar{s}=c_0<a_1
\end{array}\right.
$$
and
$$u_s\cdot a=\left\{\begin{array}{l}
sa_1\langle a_2\cdots a_na\rangle \ \ \ \ \  \ \ \ \ \ \ \  \ \ \ \ \ \  \ \ \ \ \ \ \ \ \ \ \mbox{if} \ a\geq a_1,\\
saa_1a_2\cdots a_n -a_1aa_2\cdots a_n\cdot s\ \ \ \ \ \ \ \  \mbox{if} \ a<a_1.\\
\end{array}\right.
$$
Therefore,
$$
\overline{u_s\cdot a}=\left\{\begin{array}{l}
c_0\langle c_1\cdots c_ka_1a_2\cdots a_na\rangle \ \ \ \ \ \ \mbox{if} \ \bar{s}=c_0c_1\cdots c_k,\\
c_0\langle a_1a_2\cdots a_na\rangle \ \ \ \ \ \ \ \ \ \ \ \ \ \  \ \ \mbox{if} \ \bar{s}=c_0>a_1,\\
a_1\langle c_0a_2\cdots a_n a\rangle\ \ \ \ \ \ \ \ \ \ \ \ \ \ \ \ \mbox{if} \ \bar{s}=c_0<a_1.
\end{array}\right.
$$

If $u_s=a_0a_1\cdots a_n\cdot s$, then $\overline{u_s}=a_0\langle a_1\cdots
a_n\overline{s^{(0)}}\rangle$ and $\overline{u_s\cdot a}=a_0\langle a_1\cdots a_n\overline{s^{(0)}}a\rangle$.

Since $\overline{u_s}<w$, in both cases we have $\overline{u_s\cdot a}<\overline{w\cdot a}$. \ \ \ $\Box$

\begin{definition}\label{d3.2}
Let $f$ and $g$ be momic polynomials of $\mathcal{L}_{(2)}(X)$ and
$\alpha$ and $\beta$ are the coefficients of $\overline{f^{(0)}}$
and $\overline{g^{(0)}}$ respectively. We define seven different
types of compositions as follow:

\begin{enumerate}
\item[1.] If $\bar{f}=a_0a_1\cdots a_n$, $\bar{g}=a_0b_1\cdots b_m$, $(n,m\geq0)$
and $lcm(AB)\neq \langle a_1\cdots a_nb_1\cdots b_m\rangle$, where
$lcm(AB)$ denotes the least common multiple in $[X]$ of associative
words $a_1\cdots a_n$ and $b_1\cdots b_m$, then let $w=a_0\langle
lcm(AB) \rangle$. The composition of type I of $f$ and $g$ relative
to $w$ is defined by
$$
C_I( f,g)_{w}=f\langle\frac{lcm(AB)}{a_1\cdots
a_n}\rangle-g\langle\frac{lcm(AB)}{b_1\cdots b_m}\rangle.
$$

\item[2.] If $\bar{f}=\overline{f^{(1)}}=a_0a_1\cdots a_n$, $\overline{g^{(0)}}=a_i$
for some $i\geq2$ or  $\overline{g^{(0)}}=a_1$ and $a_0>a_2$, then
let $w=\bar{f}$ and the  composition of type II of $f$ and $g$
relative to $w$ is defined by
$$
C_{II}( f,g)_{w}=f-\beta^{-1}a_0a_1\cdots \hat{a_i}\cdots a_n\cdot
g,
$$
where $a_0a_1\cdots \hat{a_i}\cdots a_n=a_0a_1\cdots
a_{i-1}a_{i+1}\cdots a_n$.

\item[3.] If $\bar{f}=\overline{f^{(1)}}=a_0a_1\cdots a_n$, $\bar{g}=\overline{g^{(0)}}
=a_1$ and $a_0\leq a_2$ or $n=1$, then let $w=\bar{f}$ and the
composition of type III of $f$ and $g$ relative to $w$ is defined by
$$
C_{III}( f,g)_{\bar{f}}=f+ga_0a_2\cdots a_n.
$$

\item[4.] If $\bar{f}=\overline{f^{(1)}}=a_0a_1\cdots a_n$, $g^{(1)}\neq0$,
$\overline{g^{(0)}}=a_1$ and $a_0\leq a_2$ or $n=1$, then for any $a<a_0$ and
 $w=a_0\langle a_1\cdots a_na\rangle$, the  composition of type IV
of $f$ and $g$ relative to $w$ is defined by
$$
C_{IV}( f,g)_{w}=fa-\beta^{-1}a_0aa_2\cdots a_n\cdot g.
$$

\item[5.] If $\bar{f}=\overline{f^{(1)}}=a_0a_1\cdots a_n$, $f^{(0)}\neq0$, $g^{(1)}
\neq0$ and $\overline{g^{(0)}}=b\notin \{a_i\}_{i=1}^n$, then let
$w=a_0\langle a_1 \cdots a_nb\rangle$ and the  composition of type V
of $f$ and $g$ relative to $w$ is defined by
$$
C_{V}( f,g)_{w}=fb-\beta^{-1}a_0a_1\cdots a_n\cdot g.
$$

\item[6.] If $\overline{f^{(0)}}=\overline{g^{(0)}}=a$ and $f^{(1)}\neq0$, then for
any $a_0a_1\in R$ and $w=a_0\langle a_1a\rangle$, the  composition
of type VI of $f$ and $g$ relative to $w$ is defined by
$$
C_{VI}( f,g)_{w}=(a_0a_1)(\alpha^{-1}f-\beta^{-1}g).
$$

\item[7.] If $f^{(1)}\neq0$, $g^{(1)}\neq0$ and $\overline{f^{(0)}}=a>\overline{g^{(0)}}=b$,
 then for any $a_0>a$ and $w=a_0ba$, the  composition of type VII
of $f$ and $g$ relative to $w$ is defined by
$$
C_{VII}( f,g)_{w}=\alpha^{-1}(a_0b)f-\beta^{-1}(a_0a)g.
$$

\end{enumerate}
\end{definition}

Immediately, we have $\overline{C_\lambda(f,g)_{w}}< w$.

\ \

\noindent {\bf Remark:} In the paper of V.V. Talapov \cite{ta}, only
the compositions of types I, II and III  are defined.

\begin{definition}
Given a set $S$ of monic polynomials of $\mathcal{L}_{(2)}(X)$ and
$w\in N$, a polynomial $f\in \mathcal{L}_{(2)}(X)$ is called trivial
modulo $S$ and $w$, denoted by $f\equiv0 \ \ mod(S,w)$, if $f$ is a
linear combination of normal $S$-words whose leading words are less
than $w$, i.e., $f=\sum_{i}\alpha_{i}u_{s_{i}}$, where $\alpha_{i}
\in {\bf k}$, $u_{s_{i}}$ are normal $S$-words and
$\overline{u_{s_{i}}}<w$. For any $f,g\in \mathcal{L}_{(2)}(X)$, we
say $f\equiv g \ \ mod(S,w)$ if $f-g\equiv0 \ \ mod(S,w)$.

The set $S$ is a Gr\"{o}bner-Shirshov basis in
$\mathcal{L}_{(2)}(X)$ if $S$ is closed under compositions, which
means every composition of any two elements of $S$ is trivial modulo
$S$ and corresponding $w$, i.e., $(\forall f,g\in S)\
C_\lambda(f,g)_w\equiv0 \ mod(S,w)$.
\end{definition}

\begin{lemma}\label{l3.6}
If $sa_1a_2\cdots a_n$ is a normal $s$-word with leading word $w$,
then for any $a_{i_1}<\bar{s}$,
$$sa_1a_2\cdots a_n\equiv sa_{i_1}a_{i_2}\cdots a_{i_n} \ \ mod(s,w),$$
where $\langle a_{i_1}a_{i_2}\cdots a_{i_n}\rangle=a_1a_2\cdots
a_n$.
\end{lemma}
{\bf Proof.} There is nothing to prove if $a_{i_1}=a_1$. Suppose
that $a_{i_1}=a_j>a_1$ for some $j\geq2$. Then we have
\begin{eqnarray*}
&&sa_1a_2\cdots a_n\\
&=&sa_1a_ja_2\cdots \hat{a_i} \cdots a_n\\
&=&sa_ja_1a_2\cdots \hat{a_i} \cdots a_n+(a_ja_1)a_2\cdots \hat{a_i}
\cdots a_n\cdot s.
\end{eqnarray*}
Since $a_{i_1}<\bar{s}$, it is easy to see that
$\overline{(a_ja_1)a_2\cdots \hat{a_i} \cdots a_n\cdot
s}<\overline{sa_1a_2\cdots a_n}=w$. The result follows.   $\Box$

The following lemma plays a key role in this paper.

\begin{lemma}\label{l3.4} Let $S$ be a Gr\"{o}bner-Shirshov basis in
$\mathcal{L}_{(2)}(X)$. If
$w=\overline{u_{s_1}}=\overline{u_{s_2}}$, where $s_1,s_2\in S$ and
$u_{s_1},u_{s_2}$ are normal $S$-words, then for some
$0\neq\alpha\in {\bf k}$,
$$
u_{s_1}\equiv \alpha u_{s_2} \ mod(S,w).
$$
\end{lemma}
{\bf Proof.} There are three main cases to consider.

Case 1. $u_{s_1}=s_1a_1a_2\cdots a_n$, $u_{s_2}=s_2b_1b_2\cdots b_m$.

(1.1) If $\bar{s}_1=\overline{s^{(1)}_1}=c_0c_1\dots c_k$ and $\bar{s}_2=
\overline{s^{(1)}_2}=d_0d_1\dots d_l$, then $c_0=d_0$ and
$$w=c_0\langle c_1\cdots c_ka_1a_2\cdots a_n\rangle=d_0\langle d_1\cdots
d_lb_1b_2\cdots b_m\rangle=c_0\langle lcm(CD)T\rangle,$$ where $T\in
[X]$ such that $\langle c_1\cdots c_ka_1a_2\cdots a_n\rangle=\langle
d_1\cdots d_lb_1b_2\cdots b_m\rangle=\langle lcm(CD)T\rangle$. Thus,
By Lemmas \ref{l3.6} and  \ref{l3} we have
\begin{eqnarray*}
&&s_1a_1a_2\cdots a_n-s_2b_1b_2\cdots b_m\\
&=&s_1\langle \frac{lcm(CD)}{c_1\cdots c_k}T\rangle-s_2\langle \frac{lcm(CD)}{d_1\cdots d_l}T\rangle\\
&\equiv&(s_1\langle \frac{lcm(CD)}{c_1\cdots c_k}\rangle-s_2\langle
\frac{lcm(CD)}{d_1\cdots d_l}\rangle)
\langle T\rangle\\
&\equiv&C_I(s_1,s_2)_{w'}\langle T\rangle\\
&\equiv& 0 \ \ mod(S,w),
\end{eqnarray*}
where $w'=c_0\langle lcm(CD)\rangle$ and $w=\overline{w'\langle
T\rangle}$.

(1.2) If $\bar{s}_1=\overline{s^{(1)}_1}=c_0c_1\dots c_k$ and $\bar{s}_2=\overline{s^{(0)}_2}=d$,
 then there are two subcases to be discussed.

(1.21) If  $d>b_1$ then
$$w=c_0\langle c_1\cdots c_ka_1a_2\cdots a_n\rangle=db_1b_2\cdots b_m,$$
which implies $c_0=d$ and $\langle c_1\cdots c_ka_1a_2\cdots a_n\rangle=b_1b_2\cdots b_m$.

Hence,
\begin{eqnarray*}
&&s_1a_1a_2\cdots a_n-s_2b_1b_2\cdots b_m\\
&\equiv&s_1a_1a_2\cdots a_n-(s_2c_1\cdots c_k)a_1a_2\cdots a_n\\
&\equiv&(s_1-s_2c_1\cdots c_k)a_1a_2\cdots a_n\\
&\equiv&C_I(s_1,s_2)_{\bar{s}_1}a_1a_2\cdots a_n\\
&\equiv& 0 \ \ mod(S,w).
\end{eqnarray*}

(1.22) If  $d<b_1$ then $a_1\geq c_1$. In fact, if $a_1< c_1\
(<c_0)$, then $w=c_0a_1\langle c_1 \cdots c_ka_2\cdots
a_n\rangle=b_1db_2\cdots b_m,$ which implies $c_0=b_1$, $a_1=d$ and
$\langle c_1\cdots c_ka_2\cdots a_n\rangle=b_2\cdots b_m$.
 This is impossible because $c_1<c_0=b_1\leq b_i \ (2\leq i\leq m)$. Thus we have $a_1\geq c_1$ and
$$w=c_0c_1\langle c_2\cdots c_ka_1a_2\cdots a_n\rangle=b_1db_2\cdots b_m,$$
which implies $c_0=b_1$, $c_1=d$ and $\langle c_2\cdots c_ka_1a_2\cdots a_n\rangle=b_2\cdots b_m$.

By noting that $c_0=b_1\leq b_i=c_2$ for some $2\leq i \leq m$, we have
\begin{eqnarray*}
&&s_1a_1a_2\cdots a_n+s_2b_1b_2\cdots b_m\\
&=&s_1a_1a_2\cdots a_n+(s_2c_0c_2\cdots c_k)a_1a_2\cdots a_n\\
&=&(s_1+s_2c_0c_2\cdots c_k)a_1a_2\cdots a_n\\
&\equiv&C_{III}(s_1,s_2)_{\bar{s}_1}a_1a_2\cdots a_n\\
&\equiv& 0 \ \ mod(S,w).
\end{eqnarray*}

(1.3) If $\bar{s}_1=\overline{s^{(0)}_1}=c$ and
$\bar{s}_2=\overline{s^{(0)}_2}=d$, then we have $n=m$. Thus, we may
 assume that $n=m\geq1$. There are two subcases to consider.

(1.31) If either $c>a_1, \ d>b_1$ or $c<a_1, \ d<b_1$, then
$$w=ca_1\cdots a_n=db_1\cdots b_m$$
or
$$w=a_1ca_2\cdots a_n=b_1db_2\cdots b_m,$$
which implies $c=d$, $a_i=b_i \ (\forall i)$ and $n=m$.

It is easy to see that
\begin{eqnarray*}
&&s_1a_1a_2\cdots a_n-s_2b_1b_2\cdots b_m\\
&=&(s_1-s_2)a_1\cdots a_n\\
&=&C_I(s_1,s_2)_{\bar{s}_1}a_1\cdots a_n\\
&\equiv& 0 \ \ mod(S,w).
\end{eqnarray*}

(1.32) If  $c>a_1$ but $d<b_1$, then
$$w=ca_1\cdots a_n=b_1db_2\cdots b_m,$$
which implies $c=b_1$, $d=a_1$, $a_i=b_i \ (i\geq2)$ and $n=m$.

Obviously,
\begin{eqnarray*}
&&s_1a_1a_2\cdots a_n+s_2b_1b_2\cdots b_m\\
&=&(s_1\bar{s}_2-\bar{s}_1s_2)a_2\cdots a_n\\
&=&(s_1(\bar{s}_2-s_2)-(\bar{s}_1-s_1)s_2)a_2\cdots a_n\\
&\equiv& 0 \ \ mod(S,w).
\end{eqnarray*}

Case 2. $u_{s_1}=s_1a_1a_2\cdots a_n$, $u_{s_2}=b_0b_1b_2\cdots b_m\cdot s_2$. We may assume that $s_2$
 is ${(0)}$-monic and $\overline{s^{(0)}_2}=d$. Then $w=b_0\langle b_1\cdots  b_md\rangle$.

(2.1) If $\bar{s}_1=\overline{s^{(1)}_1}=c_0c_1\dots c_k$, then $c_0=b_0$ and
$$w=c_0\langle c_1\cdots c_ka_1a_2\cdots a_n\rangle=b_0\langle b_1\cdots  b_md\rangle.$$

(2.11) If $d \notin \{c_i\}_{i=1}^k$, then there exists an $a_i \
(1\leq i\leq n)$ such that $d=a_i$. Thus,
\begin{eqnarray*}
&&s_1a_1a_2\cdots a_n-b_0b_1b_2\cdots b_m\cdot s_2\\
&\equiv&(s_1a_i)a_1a_2\cdots \hat{a_i} \cdots a_n-(c_0c_1\cdots c_k\cdot s_2)a_1a_2\cdots \hat{a_i} \cdots a_n\\
&\equiv&(s_1\overline{s^{(0)}_2}-\bar{s}_1s_2)a_1a_2\cdots \hat{a_i} \cdots a_n.\\
\end{eqnarray*}
If $s_2^{(1)}=0$, then
\begin{eqnarray*}
&&(s_1\overline{s^{(0)}_2}-\bar{s}_1s_2)a_1a_2\cdots \hat{a_i} \cdots a_n\\
&=&(s_1\bar{s}_2-\bar{s}_1s_2)a_1a_2\cdots \hat{a_i} \cdots a_n\\
&\equiv& 0 \ \ mod(S,w).
\end{eqnarray*}

If $s_1^{(0)}=0$, i.e., $s_1=s_1^{(1)}=\bar{s}_1+r_1^{(1)}$, then let $s^{(0)}_2=\overline{s^{(0)}_2}+r_2^{(0)}$
and we have
\begin{eqnarray*}
&&s_1\overline{s^{(0)}_2}-\bar{s}_1s_2\\
&=&(\bar{s}_1+r_1^{(1)})\overline{s^{(0)}_2}-\bar{s}_1s^{(0)}_2\\
&=&r_1^{(1)}\overline{s^{(0)}_2}-\bar{s}_1r_2^{(0)}\\
&=&r_1^{(1)}\overline{s^{(0)}_2}-\bar{s}_1r_2^{(0)}+r_1^{(1)}r_2^{(0)}-r_1^{(1)}r_2^{(0)}\\
&=&r_1^{(1)}s^{(0)}_2-s_1r_2^{(0)}\\
&=&r_1^{(1)}s_2-s_1r_2^{(0)},
\end{eqnarray*}
which implies $(s_1\overline{s^{(0)}_2}-\bar{s}_1s_2)a_1a_2\cdots \hat{a_i} \cdots a_n\equiv 0 \ \ mod(S,w)$
immediately.

If $s_2^{(1)}\neq0$ and $s_1^{(0)}\neq0$, then
\begin{eqnarray*}
&&(s_1\overline{s^{(0)}_2}-\bar{s}_1s_2)a_1a_2\cdots \hat{a_i} \cdots a_n\\
&\equiv&C_V(s_1,s_2)_{w'}a_1a_2\cdots \hat{a_i} \cdots a_n \\
&\equiv& 0 \ \ mod(S,w),
\end{eqnarray*}

where $w'=c_0\langle c_1\cdots c_kd\rangle$ and $w=\overline{w'a_1a_2\cdots \hat{a_i} \cdots a_n}$.

(2.12) If $d=c_i$ for some $i\geq2$, or $d=c_1$ and $c_0>c_2$, then
\begin{eqnarray*}
&&s_1a_1a_2\cdots a_n-b_0b_1b_2\cdots b_m\cdot s_2\\
&\equiv&s_1a_1a_2 \cdots a_n-(c_0c_1\cdots \hat{c_i}\cdots c_k\cdot s_2)a_1a_2 \cdots a_n\\
&\equiv&(s_1-c_0c_1\cdots \hat{c_i}\cdots c_k\cdot s_2)a_1a_2\cdots a_n\\
&\equiv&C_{II}(s_1,s_2)_{\bar{s}_1}a_1a_2 \cdots a_n \\
&\equiv& 0 \ \ mod(S,w),
\end{eqnarray*}
where $c_i=d$.

(2.13) If $d=c_1$ and $c_0\leq c_2$, then by the form of $w$, we have
$b_0b_1\cdots b_m=c_0\langle c_2\cdots c_ka_1\cdots a_n\rangle\in R$, which implies $c_2\geq c_0>a_1$. Thus,
\begin{eqnarray*}
&&s_1a_1a_2\cdots a_n-b_0b_1b_2\cdots b_m\cdot s_2\\
&=&s_1a_1a_2\cdots a_n-c_0a_1\langle c_2\cdots c_ka_2\cdots a_n\rangle\cdot s_2\\
&=&(s_1a_1-c_0a_1c_2\cdots c_k\cdot s_2)a_2\cdots a_n\\
&=&C_{IV}(s_1,s_2)_{w'}a_2\cdots a_n \\
&\equiv& 0 \ \ mod(S,w),
\end{eqnarray*}
where $w'=c_0\langle c_1\cdots c_ka_1\rangle$ and $w=w'a_2\cdots a_n$.

(2.2) If $\bar{s}_1=\overline{s^{(0)}_1}=c$ and $\overline{s^{(0)}_2}=d$, then $n=m+1\geq2$ since
$w=b_0\langle b_1\cdots  b_md\rangle$ and $m\geq 1$.

(2.21) If $c>a_1$, then $w=ca_1\cdots a_n=b_0\langle b_1\cdots  b_md\rangle$, which implies $b_0=c$.

(2.211) If $d\geq b_1$, then $a_1=b_1$, $a_2\cdots a_n=\langle b_2\cdots b_md\rangle$ and
\begin{eqnarray*}
&&s_1a_1a_2\cdots a_n-b_0b_1b_2\cdots b_m\cdot s_2\\
&=&(s_1b_1d)b_2\cdots b_m-((b_0b_1)\cdot s_2)b_2\cdots b_m\\
&=&(s_1b_1d-(b_0b_1)\cdot s_2)b_2\cdots b_m\\
&=&(s_1b_1\overline{s^{(0)}_2}-(\bar{s}_1b_1)\cdot s_2^{(0)})b_2\cdots b_m\\
&=&(s_1b_1(\overline{s^{(0)}_2}-s_2^{(0)})-((\bar{s}_1-s_1)b_1)\cdot s_2^{(0)})b_2\cdots b_m\\
&=&(s_1b_1)\langle r_2^{(0)}b_2\cdots b_m\rangle-(r_1b_1)b_2\cdots b_m\cdot s_2 \\
&\equiv& 0 \ \ mod(S,w),
\end{eqnarray*}
where $s_2^{(0)}=\overline{s^{(0)}_2}+r_2^{(0)}$ and $s_1=\bar{s}_1+r_1$.

(2.212) If $d<b_1$, then $w=cdb_1\cdots b_m$. Suppose that
$s_1=c+\sum_{c_i<c}\alpha_i c_i$, $s_2^{(0)}=d+\sum_{d_j<d}\beta_j
d_j$. Thus,
\begin{eqnarray*}
&&s_1a_1a_2\cdots a_n-b_0b_1b_2\cdots b_m\cdot s_2\\
&=&s_1db_1b_2\cdots b_m-cb_1b_2\cdots b_m\cdot s_2\\
&=&(s_1db_1-(cb_1)\cdot s_2)b_2\cdots b_m.\\
&=&(s_1db_1-(s_1b_1) s_2+\sum_{c_i<c} \alpha_i(c_ib_1)s_2)b_2\cdots b_m\\
&=&(s_1b_1d+(b_1d)s_1-(s_1b_1)s_2+\sum_{c_i<c}  \alpha_i(c_ib_1)s_2)b_2\cdots b_m\\
&=&(s_1b_1(d-s_2)+(b_1d)s_1+\sum_{c_i<c} \alpha_i(c_ib_1) s_2)b_2\cdots b_m\\
&=&(-\sum_{d_j<d} \beta_j s_1b_1d_j+(b_1d)s_1+(\sum_{c_i<c}\alpha_i c_ib_1)\cdot s_2)b_2\cdots b_m\\
&=&(-\sum_{d_j<d} \beta_j s_1d_jb_1+\sum_{d_j<d}\beta_j(b_1d_j)s_1+(b_1d)s_1+(\sum_{c_i<c}\alpha_i c_ib_1)
\cdot s_2)b_2\cdots b_m\\
&\equiv& 0 \ \ mod(S,w).
\end{eqnarray*}

(2.22) If $c<a_1$, then $w=a_1ca_2\cdots a_n=b_0\langle b_1\cdots  b_md\rangle$ and $a_1=b_0$.
In this case, $d\geq b_1$, and then $b_1=c$, $d=a_i$ for some $i\geq2$. Otherwise, if $d<b_1$,
then $d=c$. This implies $a_i=b_{i-1}$ for any $i\geq1$ and $b_0=a_1\leq a_2=b_1$, which is a
contradiction. Therefore,
\begin{eqnarray*}
&&s_1a_1a_2\cdots a_n+b_0b_1b_2\cdots b_m\cdot s_2\\
&=&-(a_1s_1)a_2\cdots a_n+a_1b_1b_2\cdots b_m\cdot s_2\\
&=&-((a_1s_1)d)a_2\cdots \hat{a_i}\cdots a_n+(a_1c)a_2\cdots \hat{a_i}\cdots a_n\cdot s_2\\
&=&((s_1a_1)\overline{s^{(0)}_2}+(a_1\bar{s}_1)\cdot s_2)a_2\cdots \hat{a_i}\cdots a_n\\
&=&((s_1a_1)(\overline{s^{(0)}_2}-s_2^{(0)})+(a_1(\bar{s}_1-s_1))\cdot s_2)a_2\cdots \hat{a_i}\cdots a_n\\
&\equiv& 0 \ \ mod(S,w).
\end{eqnarray*}

Case 3. $u_{s_1}=a_0a_1a_2\cdots a_n\cdot s_1$,
$u_{s_2}=b_0b_1b_2\cdots b_n\cdot s_2$. We may assume that both
$s_1$ and $s_2$ are ${(0)}$-monic. Suppose that
$\overline{s^{(0)}_1}=c$ and $\overline{s^{(0)}_2}=d$. Then
$w=a_0\langle a_1a_2\cdots a_nc\rangle=b_0\langle b_1b_2\cdots
b_nd\rangle$ and $a_0=b_0$.

(3.1) If $c=d$, then $a_i=b_i$ for all $i$ and
$$
a_0a_1a_2\cdots a_n\cdot s_1-b_0b_1b_2\cdots b_n\cdot s_2=a_0a_1a_2\cdots a_n\cdot (s_1-s_2).
$$

If $s^{(1)}_1=s^{(1)}_2=0$, i.e., $\bar{s}_1=\overline{s^{(0)}_1}=\overline{s^{(0)}_2}=\bar{s}_2=c$, then
$$
a_0a_1a_2\cdots a_n\cdot (s_1-s_2)=a_0a_1a_2\cdots a_n\cdot C_I(s_1,s_2)\equiv 0. \ \ \ \ mod(S,w).
$$

If $s^{(1)}_1\neq0$, then
\begin{eqnarray*}
&&a_0a_1a_2\cdots a_n\cdot (s_1-s_2)\\
&=&((a_0a_1)(s_1-s_2))a_2\cdots a_n \\
&=&C_{VI}(s_1,s_2)_{w'}a_2\cdots a_n\\
&\equiv& 0 \ \ mod(S,w),
\end{eqnarray*}
where $w'=a_0\langle a_1c\rangle$.

(3.2) If $c\neq d$, say, $c>d$, then $w=a_0\langle cda_1\cdots \hat{a_i}\cdots
a_n\rangle=a_0\langle cdb_1\cdots \hat{b_j}\cdots b_n\rangle$ for some $a_i$ and $b_j$.

(3.21) If $d\geq b_1$, then $w=a_0b_1\langle cdb_2\cdots \hat{b_j}\cdots b_n
\rangle=a_0a_1\langle cda_2\cdots \hat{a_i}\cdots a_n\rangle$, which implies
 $a_1=b_1$, $a_2\cdots \hat{a_i}\cdots a_n=b_2\cdots \hat{b_j}\cdots b_n$. Thus,
\begin{eqnarray*}
&&a_0a_1a_2\cdots a_n\cdot s_1-b_0b_1b_2\cdots b_n\cdot s_2\\
&=&((a_0b_1d)\cdot s_1)a_2\cdots \hat{a_i}\cdots a_n-((a_0b_1c)\cdot s_2)b_2\cdots\hat{b_j}\cdots b_n \\
&=&(a_0b_1d\cdot s_1-a_0b_1c\cdot s_2)b_2 \cdots\hat{b_j}\cdots b_n \\
&=&(a_0b_1(d-s_2)\cdot s_1-a_0b_1(c-s_1)\cdot s_2)b_2 \cdots\hat{b_j}\cdots b_n \\
&\equiv& 0 \ \ mod(S,w).
\end{eqnarray*}

(3.22) If $d<b_1$, then $w=a_0db_1\cdots b_n=a_0a_1\langle a_2\cdots a_nc\rangle$,
which implies $a_1=d$ and $c=b_i$ for some $i$.

(3.221) If $c=b_1<a_0$, then $a_i=b_i \ (i\geq2)$ and $w=a_0dcb_2\cdots b_n$. We have
\begin{eqnarray*}
&&a_0a_1a_2\cdots a_n\cdot s_1-b_0b_1b_2\cdots b_n\cdot s_2\\
&=&((a_0d)\cdot s_1)a_2\cdots a_n-((a_0c)\cdot s_2)a_2\cdots a_n \\
&=&(a_0d\cdot s_1-a_0c\cdot s_2)a_2 \cdots a_n.
\end{eqnarray*}
If $s^{(1)}_1=0$, then we may suppose that
$s_1=c+\sum_{c_i<c}\alpha_ic_i$ and
$s_2^{(0)}=d+\sum_{d_j<d}\beta_jd_j$. We have
\begin{eqnarray*}
&&(a_0d s_1-a_0cs_2)a_2 \cdots a_n \\
&=&((a_0s_1)d+s_1da_0-a_0c\cdot s_2)a_2 \cdots a_n \\
&=&((a_0s_1)s_2-a_0c\cdot s_2+s_1da_0+\sum_{d_j<d}\beta_js_1a_0d_j)a_2 \cdots a_n\\
&=&((a_0(s_1-c)s_2+s_1da_0+\sum_{d_j<d}\beta_js_1d_ja_0d_j-\sum_{d_j<d}\beta_j(a_0d_j)s_1)a_2 \cdots a_n\\
&=&(\sum_{c_i<c}\alpha_i(a_0c_i)s_2+s_1da_0+\sum_{d_j<d}\beta_js_1d_ja_0d_j-\sum_{d_j<d}\beta_j(a_0d_j)s_1)a_2 \cdots a_n\\
&\equiv& 0 \ \ mod(S,w).
\end{eqnarray*}
If $s^{(1)}_2=0$, then we have
\begin{eqnarray*}
&&(a_0d s_1-a_0c\cdot s_2)a_2 \cdots a_n \\
&=&(a_0d s_1-a_0s_2c-s_2ca_0)a_2 \cdots a_n \\
&=&(a_0(d-s_2) s_1-a_0s_2(c-s_1)-s_2ca_0)a_2 \cdots a_n \\
&\equiv& 0 \ \ mod(S,w).
\end{eqnarray*}
If $s^{(1)}_i\neq0 \ (i=1,2)$, then let $w'=a_0dc$. We have
$w=w'a_2\cdots a_n$ and
\begin{eqnarray*}
&&(a_0d s_1-a_0c\cdot s_2)a_2 \cdots a_n \\
&=&C_{VII}(s_2,s_1)_{w'}a_2 \cdots a_n \\
&\equiv& 0 \ \ mod(S,w).
\end{eqnarray*}

(3.222) If $c=b_i>b_1$ for some $i\geq2$, then
\begin{eqnarray*}
&&a_0a_1a_2\cdots a_n\cdot s_1-b_0b_1b_2\cdots b_n\cdot s_2\\
&=&((a_0db_1)\cdot s_1-(a_0b_1c)\cdot s_2)b_2\cdots\hat{b_j}\cdots b_n \\
&=&(a_0b_1d\cdot s_1+b_1da_0\cdot s_1-a_0b_1c\cdot s_2)b_2\cdots\hat{b_j}\cdots b_n \\
&\equiv&(a_0b_1\overline{s^{(0)}_2}\cdot s_1-a_0b_1\overline{s^{(0)}_1}\cdot s_2)b_2\cdots\hat{b_j}\cdots b_n \\
&\equiv&(((a_0b_1)\cdot s_2)\cdot s_1-((a_0b_1)\cdot s_1)\cdot s_2)b_2\cdots\hat{b_j}\cdots b_n \\
&\equiv& 0 \ \ mod(S,w).
\end{eqnarray*}

The proof is complete.\ \ \ $\Box$

\begin{theorem}{\bf (Composition-Diamond lemma for metabelian Lie algebras)} \label{CDML} Let
$S\subset{ \mathcal{L}_{(2)}(X)}$ be a nonempty set of monic
polynomials and $Id(S)$ be the ideal of $\mathcal{L}_{(2)}(X)$
generated by $S$. Then the following statements are equivalent.
\begin{enumerate}
\item[(i)] $S$ is a Gr\"{o}bner-Shirshov basis.
\item[(ii)] $f\in{Id(S)}\Rightarrow\bar{f}=\overline{u_s}$ for
some normal $S$-word $u_s$.
\item[(iii)]$Irr(S)=\{u \ | \ u \in N, u\neq\overline{v_s} \ \mbox{for any normal S-word}  \ v_s\}$
is a $\bf k$-basis for
$\mathcal{L}_{(2)}(X|S)=\mathcal{L}_{(2)}(X)/Id(S)$.
\end{enumerate}
\end{theorem}

{\bf Proof.} $(i)\Rightarrow (ii)$. Let $S$ be a
Gr\"{o}bner-Shirshov basis and $0\neq f\in Id(S).$ Then
 by Lemma \ref{l3.2} $f$ has an expression
 $f=\sum\alpha_iu_{s_i}$, where $0\neq\alpha_i\in {\bf k}, \
 u_{s_i}$ are normal $S$-words. Denote
 $w_i=\overline{u_{s_i}}, \ i=1,2,\dots$.
 We may assume without loss of generality that
$$w_1=w_2=\cdots =w_l> w_{l+1}\geq w_{l+2}\geq \cdots$$
\noindent for some $l\geq 1$.

The claim of the theorem is obvious if $l=1$.

Now suppose that $l>1$. Then
$\overline{u_{s_1}}=w_1=w_2=\overline{u_{s_2}}$. By Lemma \ref{l3.4}, for some $\alpha\in {\bf k}$,
$$
u_{s_2}\equiv \alpha u_{s_1} \ mod(S,w_1).
$$
Thus,
\begin{eqnarray*}
&&\alpha_1u_{s_1}+\alpha_2u_{s_2}\\
&=&(\alpha_1+\alpha\alpha_2)u_{s_1}+\alpha_2(u_{s_2}-\alpha u_{s_1})\\
&\equiv&(\alpha_1+\alpha\alpha_2)u_{s_1} \ \ \ \ mod(S,w_1).
\end{eqnarray*}
Therefore, if $\alpha_1+\alpha\alpha_2\neq 0$ or $l>2$, then the
result follows from the induction on $l$. For the case
$\alpha_1+\alpha\alpha_2= 0$ and $l=2$, we use the induction on
$w_1$. Now the result follows.

$(ii)\Rightarrow (iii).$ For any $f\in \mathcal{L}_{(2)}(X)$, we have
\begin{equation*}
f=\sum\limits_{\overline{ u_{s_i}}\leq\bar f}\alpha_i u_{s_i}+
\sum\limits_{\overline{v_j}\leq \bar f}\beta_jv_j,
\end{equation*}
where  $\alpha_i,\beta_j\in {\bf k},   \ v_j \in Irr(S)$ and
$u_{s_i}$ are normal $S$-words. Therefore, the set $Irr(S)$
generates the algebra $\mathcal{L}_{(2)}(X)/Id(S)$.

On the other hand, suppose that $h=\sum\alpha_iv_i=0$ in
$\mathcal{L}_{(2)}(X)/Id(S)$, where $\alpha_i\in {\bf k}$, $v_i\in
{Irr(S)}$. This means that $h\in{Id(S)}$. Then all $\alpha_i$ must
be equal to zero. Otherwise, $\overline{h}=v_j$ for some $j$ which
contradicts $(ii)$.
 \\

$(iii)\Rightarrow (i).$
 For any $f,g\in{S}$,  we
have
$$
C_\lambda(f,g)_{w}=\sum\limits_{\overline{ u_{s_i}}< w}\alpha_i
u_{s_i}+ \sum\limits_{\overline{v_j}< w}\beta_jv_j.
$$
Since $C_\lambda(f,g)_{w}\in {Id(S)}$ and by
$(iii)$, we have
$$
C_\lambda(f,g)_{w}=\sum\limits_{\overline{ u_{s_i}}< w}\alpha_i u_{s_i}.
$$
Therefore, $S$ is a Gr\"{o}bner-Shirshov basis. \ \ \ \ $\Box$

\begin{lemma} \label{l4.1}(\cite{ta}) Suppose that $f\in \mathcal{L}_{(2)}(X)$. Then
there exists an element $f'\in \mathcal{L}_{(2)}(X)$ such that
$Id(f)=Id(f')$, $\bar{f'}\leq\bar{f}$, $f'^{(0)}=f^{(0)}$ and no
word occurring in $f'^{(1)}$ contains
 $\overline{f^{(0)}}$ as a strict subword.
\end{lemma}
{\bf Proof.} If no word occurring in $f^{(1)}$ contains
 $\overline{f^{(0)}}$ as a strict subword, then we are done.
If $\overline{f^{(1)}}$ contains $\overline{f^{(0)}}$ as a strict subword, say
$\bar{f}=\overline{f^{(1)}}=a_0a_1\cdots a_n$, $\overline{f^{(0)}}=a_i$
for some $i\geq2$ or  $\overline{f^{(0)}}=a_1$ and $a_0>a_2$, then
let $f_1$ be the composition of type II of $f$ and itself:
$$
f_1=C_{II}( f,f)_{\bar{f}}=f-\beta^{-1}a_0a_1\cdots \hat{a_i}\cdots a_n\cdot
f,
$$
where $a_i=\overline{f^{(0)}}$.  It is obvious that $Id(f)=Id(f_1)$, and $\bar{f_1}<\bar{f}$, $f_1^{(0)}=f^{(0)}$.
If $\overline{f_1^{(1)}}$ contains $\overline{f^{(0)}}$ as a strict subword, we again consider the composition $f_2=C_{II}(f_1,f_1)_{\bar{f_1}}$, and so on. By induction on the leading word, we obtain an element $f'$ such that $Id(f)=Id(f')$, $\bar{f'}\leq\bar{f}$, $f'^{(0)}=f^{(0)}$, and either $f'=f'^{(0)}$ or $\overline{f'^{(1)}}$ dose not contain $\overline{f^{(0)}}$ as a strict subword.

Arguments analogous to the one given above for the leading word also apply to other regular $R$-words occurring in the expansion of $f$ and containing $\overline{f^{(0)}}$ as a strict subword. Finally, we have the one we want.  $\Box$

\begin{lemma} \label{l3.3} Suppose that $\bar{f}=\overline{f^{(1)}}=a_0a_1\cdots
a_n$, $g^{(1)}\neq0$, $\overline{g^{(0)}}=a_1$ and $a_0\leq a_2$ or
$n=1$. If $f^{(0)}=0$, then for $a=a_1<a_0$ and $w=a_0\langle
a_1\cdots a_na\rangle$, the composition of type IV of $f$ and $g$ is
trivial.
\end{lemma}
{\bf Proof.} We may suppose that $g$ is $(0)$-monic. Then
\begin{eqnarray*}
C_{IV}( f,g)_{w}&=&fa_1-\bar{f}\cdot g\\
&=&r_f^{(1)}\cdot\overline{g^{(0)}}-\bar{f}\cdot r_g^{(0)}\\
&=&r_f^{(1)}\cdot\overline{g^{(0)}}-\bar{f}\cdot r_g^{(0)}+r_f^{(1)}\cdot r_g^{(0)}-r_f^{(1)}\cdot r_g^{(0)}\\
&=&r_f^{(1)}(\overline{g^{(0)}}+r_g^{(0)})-(\bar{f}+r_f^{(1)})\cdot r_g^{(0)}\\
&=&r_f^{(1)}\cdot g-f\cdot r_g^{(0)}\\
&\equiv&0 \ \ \ mod(\{f,g\}, w),
\end{eqnarray*}
where $f=f^{(1)}=\bar{f}+r_f^{(1)}$ and $g^{(0)}=
\overline{g^{(0)}}+r_g^{(0)}$.     $\Box$

\begin{lemma} \label{l3.5} The compositions of type I, V and VI formed by $f$ itself are always trivial.
\end{lemma}
{\bf Proof.} For type I and VI, the result is obvious. We only check type V.
 Suppose that $\bar{f}=\overline{f^{(1)}}=a_0a_1\cdots a_n$,
 $\overline{f^{(0)}}=b\notin \{a_i\}_{i=1}^n$, and $w=a_0\langle a_1\cdots a_nb\rangle$. We have
\begin{eqnarray*}
C_{V}( f,f)_{w}&=&fb-\beta^{-1}a_0a_1\cdots a_n\cdot f\\
&=&f\cdot\overline{f^{(0)}}-\beta^{-1}\bar{f}\cdot f\\
&=&f\cdot\overline{f^{(0)}}-f\cdot\beta^{-1}(r^{(1)}+\beta \overline{f^{(0)}}+r^{(0)})\\
&=&-\beta^{-1}f\cdot(r^{(1)}+r^{(0)})\\
&=&\beta^{-1}r^{(1)}\cdot f-\beta^{-1}f\cdot r^{(0)}\\
&\equiv&0 \ \ \ mod(f, w),
\end{eqnarray*}
where $f^{(1)}=\bar{f}+r^{(1)}$ and $f^{(0)}=\beta \overline{f^{(0)}}+r^{(0)}$, $\beta\in {\bf k}$.  $\Box$

\ \

\noindent {\bf Remark:} If a subset $S$ of $\mathcal{L}_{(2)}(X)$ is
not a Gr\"{o}bner-Shirshov basis, then one can add all nontrivial
compositions of polynomials of $S$ to $S$. Continuing this process
repeatedly, we finally obtain a Gr\"{o}bner-Shirshov basis $S^C$
that generates the same ideal as $S$. Such a process is called
Shirshov's algorithm and $S^C$ is called a Gr\"{o}bner-Shirshov
complement of $S$. By Lemma \ref{l4.1}, we may assume any element of
the original relation set $S$ has no composition of type II formed
by itself and the Shirshov's algorithm do not involve compositions
discussed in Lemmas \ref{l3.3} and \ref{l3.5}.

\section{Applications}

Suppose that $\mathcal{A}$ is a metabelian Lie algebra and
$Y=\{a_i,i\in I\}\cup\{b_j,j\in J\}$ is a ${\bf k}$-basis of
$\mathcal{A}$, where $\{a_i\}$ is a basis of $\mathcal{A}^{(1)}$ and
$b_j$'s are linear independent modulo $\mathcal{A}^{(1)}$. Suppose
that $I$ and $J$ are well-ordered sets. The set of multiplications
of $Y$, say $M$, consists of the following:
\begin{eqnarray*}
&&m_{1ij}: \ a_ib_j-\sum\gamma_{ij}^ka_k,\\
&&m_{2ij}: \ b_ib_j-\sum\delta_{ij}^ka_k, \ (i>j), \\
&&m_{3ij}: \ a_ia_j, \ (i>j),
\end{eqnarray*}
where $\gamma_{ij}^k, \delta_{ij}^k\in {\bf k}$. Then we have
$\mathcal{A}=\mathcal{L}_{(2)}(Y|M)$ and since $Irr(M)=Y$, by
Theorem \ref{CDML}, $M$ is a Gr\"{o}bner-Shirshov basis for
$\mathcal{A}$ with respect to $a_i>b_j$.

Let $\mathcal{S}$ denote the free metabelian Lie product of
$\mathcal{A}$ and a
 free metabelian Lie algebra generated by a well-ordered set $X=\{x_h| h\in H\}$, i.e.,
$$\mathcal{S}=\mathcal{A}\ast \mathcal{L}_{(2)}(X)=\mathcal{L}_{(2)}(X\cup Y|M).$$

\begin{theorem} \label{t4.1} Let the notion be as above. Then with
respect to $x_h>a_i>b_j$, a Gr\"{o}bner-Shirshov complement $M^C$ of
$M$ in $\mathcal{L}_{(2)}(X\cup Y)$ consists of $M$ and some
$X$-homogenous polynomials without $(0)$-part, whose leading words
 are of the form $xy\cdots$ with an $a_i$ as a
strict subword, $x\in X$, $a_i,y\in Y$.
\end{theorem}
{\bf Proof.} For convenience, we call the $X$-homogenous polynomials
described in the theorem to satisfy property $P_X$.

 Since $M$ is a
Gr\"{o}bner-Shirshov basis in $\mathcal{L}_{(2)}(Y)$,
 we need to check the compositions  which are formed by $M$ itself
 and involve some elements in $X$. The possible types are VI and VII.

First, we check type VI. Suppose that
$\overline{m_{1ij}^{(0)}}=\overline{m_{1st}^{(0)}}=a_l$ and the
corresponding $w$ is of the forms $xx'a_l$, $xba_l$ and $x\langle
aa_l\rangle$ for some $x,x'\in X, \ b\in \{b_j\}$ and $a\in
\{a_i\}$.

If $w=xx'a_l$, then
\begin{eqnarray*}
C_{VI}(m_{1ij},m_{1st})_w &=& (xx')((\gamma_{ij}^l)^{-1}m_{1ij}-(\gamma_{st}^l)^{-1}m_{1st}) \\
&=&-\sum_{k<l}(\gamma_{ij}^l)^{-1}\gamma_{ij}^k xx'a_k+\sum_{k<l}(\gamma_{st}^l)^{-1}\gamma_{st}^k xx'a_k\\
&=&-\sum_{k<l}(\gamma_{ij}^l)^{-1}\gamma_{ij}^k xa_kx'+\sum_{k<l}(\gamma_{st}^l)^{-1}\gamma_{st}^k xa_kx'\\
&&+\sum_{k<l}(\gamma_{ij}^l)^{-1}\gamma_{ij}^k x'a_kx-\sum_{k<l}(\gamma_{st}^l)^{-1}\gamma_{st}^k x'a_kx\\
\end{eqnarray*}
and obviously it satisfies $P_X$.

If $w=xba_l$, then
\begin{eqnarray*}
C_{VI}(m_{1ij},m_{1st})_w &=& (xb)((\gamma_{ij}^l)^{-1}m_{1ij}-(\gamma_{st}^l)^{-1}m_{1st}) \\
&=&-\sum_{k<l}(\gamma_{ij}^l)^{-1}\gamma_{ij}^k xba_k+\sum_{k<l}(\gamma_{st}^l)^{-1}\gamma_{st}^k xba_k\\
\end{eqnarray*}
and still satisfies $P_X$.

If $w=xaa_l$, then
\begin{eqnarray*}
C_{VI}(m_{1ij},m_{1st})_w &=& (xa)((\gamma_{ij}^l)^{-1}m_{1ij}-(\gamma_{st}^l)^{-1}m_{1st}) \\
&=&-\sum_{k<l}(\gamma_{ij}^l)^{-1}\gamma_{ij}^k xaa_k+\sum_{k<l}(\gamma_{st}^l)^{-1}\gamma_{st}^k xaa_k\\
&\equiv&-\sum_{a_k<a}(\gamma_{ij}^l)^{-1}\gamma_{ij}^k xa_ka+\sum_{a_k<a}(\gamma_{st}^l)^{-1}\gamma_{st}^k xa_ka\\
&&-\sum_{a_k\geq a}(\gamma_{ij}^l)^{-1}\gamma_{ij}^k
xaa_k+\sum_{a_k\geq a}(\gamma_{st}^l)^{-1}\gamma_{st}^k xaa_k \ \ \
mod(M,w),
\end{eqnarray*}
and again the remainder satisfies $P_X$.

$C_{VI}(m_{1ij},m_{2st})_w$, $C_{VI}(m_{2ij},m_{2st})_w$ are similar
to $C_{VI}(m_{1ij},m_{1st})_w $.

Second, we check type VII. Suppose that $\overline{m_{1ij}^{(0)}}=a_p>a_q=\overline{m_{1st}^{(0)}}$ and $w=xa_qa_p$. Then
\begin{eqnarray*}
&&C_{VII}(m_{1ij},m_{1st})_w\\
 &=& (\gamma_{ij}^p)^{-1}(xa_q)m_{1ij}-(\gamma_{st}^q)^{-1}(xa_p)m_{1st}\\
&=&-\sum_{k<p}(\gamma_{ij}^p)^{-1}\gamma_{ij}^k xa_qa_k+\sum_{k<q}(\gamma_{st}^q)^{-1}\gamma_{st}^k xa_pa_k-x(a_pa_q)\\
&=&-\sum_{q\leq k<l}(\gamma_{ij}^p)^{-1}\gamma_{ij}^k xa_qa_k-\sum_{k<q}(\gamma_{ij}^p)^{-1}\gamma_{ij}^k xa_ka_q
   -\sum_{q\leq k<l}(\gamma_{ij}^p)^{-1}\gamma_{ij}^k x(a_qa_k)\\
&&+\sum_{k<q}(\gamma_{st}^l)^{-1}\gamma_{st}^k xa_ka_p+\sum_{k<q}(\gamma_{st}^l)^{-1}\gamma_{st}^k x(a_pa_k)-x(a_pa_q)\\
&\equiv&-\sum_{q\leq k<l}(\gamma_{ij}^p)^{-1}\gamma_{ij}^k
xa_qa_k-\sum_{k<q}(\gamma_{ij}^p)^{-1}\gamma_{ij}^k
+xa_ka_q\sum_{k<q}(\gamma_{st}^l)^{-1}\gamma_{st}^k xa_ka_p  \ \ \ \
mod(M,w),
\end{eqnarray*}
and the remainder has property $P_X$. One may check that
$C_{VII}(m_{1ij},m_{2st})_w$ and $C_{VII}(m_{2ij},m_{2st})_w$ are
the same as $C_{VII}(m_{1ij},m_{1st})_w$, which have property $P_X$.

Observing from above and the definition of compositions, we know
that the non-trivial compositions of polynomials satisfy $P_X$
themselves are only of type I and the results again satisfy $P_X$.
Also by the definition of compositions and property $P_X$, the
compositions of $M$ and polynomials satisfying $P_X$ are only of
type II and the results still satisfy $P_X$. The theorem is proved.
\ \ $\Box$

Observing from the proof of the above theorem, we have the following
proposition.
\begin{proposition} \label{p4.2} Let $\mathcal{A}_i=\mathcal{L}_{(2)}(X_i|S_i
)$, where $S_i\subset \mathcal{L}_{(2)}(X_i)^{(1)} $, $i=1,2$. Then
$S_1^C\cup S_2^C$ is a Gr\"{o}bner-Shirshov basis for the free
metabelian Lie product $\mathcal{A}_1\ast \mathcal{A}_2$, where
$S_i^C$ is a Gr\"{o}bner-Shirshov complement of $S_i$ in
$\mathcal{L}_{(2)}(X_i),\ i=1,2$.
\end{proposition}

\ \

Now, we consider partial commutative metabelian Lie algebras related
to some graphs.

Let $\Gamma=(V,E)$ be a graph, where $V$ is the set of vertices and
$E$ the set of edges. For $e\in E$ we call $o(e)$ the origin of $e$
and $t(e)$ the terminus. We say a metabelian Lie algebra is partial
commutative related to a graph $\Gamma=(V,E)$, denoted by
$\mathcal{ML}_{\Gamma}$, if
$$
\mathcal{ML}_{\Gamma}=\mathcal{L}_{(2)}(V  | \  [o(e), t(e)]=0, e\in E).
$$

In this section, we find Gr\"{o}bner-Shirshov bases for partial
commutative metabelian Lie algebras related to any circuits, trees
and 3-cube.

The following algorithm  gives a Gr\"{o}bner-Shirshov basis for
partial commutative metabelian Lie algebras with a finite relation
set.

\begin{algorithm} \label{a1}

Input: relations $f_1, \cdots, f_s$ of $\mathcal{L}_{(2)}(X)$, $f_i=xx'$, $F=\{f_1, \cdots, f_s\}$.

Output: a Gr\"{o}bner-Shirshov basis $H=\{h_1, \cdots, h_t\}$ for $\mathcal{L}_{(2)}(X|F)$.

Initialization: $H:=F$

While: \ \ \ \ \ \ \ \  $f_i=x_{i_0}x_{i_1}\cdots x_{i_n}, \ f_i=x_{j_0}x_{j_1}\cdots x_{j_m}$,
and $x_{i_0}=x_{j_0}$, $x_{i_1}\neq x_{j_1}$

Then Do:
 \ \ \ \ $h:=max\{x_{i_1},x_{j_1}\}min\{x_{i_1},x_{j_1}\}\langle x_{t_1}x_{t_2}\cdots x_{t_l}\rangle$

 \ \ \ \ \ \  \ \ \ \ \ \ where $\{x_{t_1},x_{t_2},\cdots, x_{t_l}\}=\{x_{i_0},x_{i_2},\cdots,
 x_{i_n}\}\cup\{x_{j_2},\cdots,x_{j_m}\}$

 \ \ \ If:  \ \ \ \ \ \ \ \ \ there is no $f_j\in H$ such that $f_j$ is a subword of $h$

 \ \ \ Do: \ \ \ \ \ \ \  $H:=H\cup \{h\}$

End

\end{algorithm}

\begin{definition} Let $n$ be a positive integer. A \emph{circuit (of length $n$)},
denoted by $Circ_n$, is a graph which the set of vertices is $\mathbf{Z}/n\mathbf{Z}$
 and the orientation is given by $n$ edges $e_{i,i+1}, \ i\in \mathbf{Z}/n\mathbf{Z}$,
 with $o(e_{i,i+1})=i$ and $t(e_{i,i+1})=i+1$.
\end{definition}

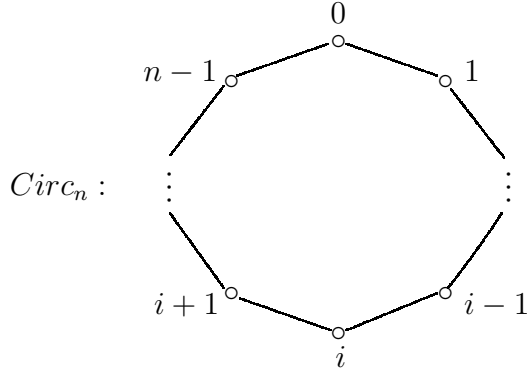
\begin{picture}(100,130)
\put(80,55){$Circ_n:$}
\put(130,98){$n-1$}\put(200,120){$0$}\put(250,98){$1$}
\put(134,10){$i+1$}\put(202,-10){$i$}\put(250,10){$i-1$}
\put(160,95){$\circ$}\put(200,110){$\circ$}\put(240,95){$\circ$}
\put(265,60){$\cdot$}\put(265,55){$\cdot$}\put(265,50){$\cdot$}
\put(138,60){$\cdot$}\put(138,55){$\cdot$}\put(138,50){$\cdot$}
\put(160,15){$\circ$}\put(200,0){$\circ$}\put(240,15){$\circ$}
\qbezier(206,112)(223,106)(240,100)
\qbezier(245,95)(255,82)(265,69)
\qbezier(244,20)(255,34)(264,48)
\qbezier(206,4)(223,11)(240,18)
\qbezier(200,4)(183,10)(166,16)
\qbezier(160,20)(150,34)(140,48)
\qbezier(140,70)(150,83)(160,96)
\qbezier(165,100)(182,106)(199,112)
\end{picture}

 \begin{theorem} For the partial commutative metabelian Lie algebra related to $Circ_n$
 $$
 \mathcal{ML}_{Circ_n}=\mathcal{L}_{(2)}(\mathbf{Z}/n\mathbf{Z}  \  | \  [i+1,
 i]=0,  \ i\in\mathbf{Z}/n\mathbf{Z}),
 $$
with the usual ordering on natural numbers, a Gr\"{o}bner-Shirshov
basis for $\mathcal{ML}_{Circ_n}$ consists of the following
relations:
\begin{eqnarray*}
&&f_{0}: \ [n-1,0]=0,\\
&&f_{i}: \ [i,i-1]=0, \  1\leq i\leq n-1, \\
&&g_{j}: \ [j,0,j+1,j+2,\cdots,n-1]=0, \ 2\leq j\leq n-2,
\end{eqnarray*}
where the brackets $[\cdots]$ is the left-normed brackets.
\end{theorem}
{\bf Proof.} The only possible compositions are of type $I$ by
$f_{n-1}, f_0$ and $g_j,f_j$, where the corresponding $w's$ are
$[n-1,0,n-2]$ and $[j,0,j-1,j+1,j+2,\cdots,n-1]$ respectively.

For the first one, $w=[n-1,0,n-2]$ and
\begin{eqnarray*}
&&C_I(f_{n-1}, f_0)_{w}\\
&=&[n-1,n-2]\cdot 0-[n-1,0,n-2]\\
&=&[n-2,0,n-1]\\
&\equiv&0 \ \ \ \ \ mod(g_{n-2},w ).
\end{eqnarray*}
For the second one, $w=[j,0,j-1,j+1,j+2,\cdots,n-1]$ and
\begin{eqnarray*}
&&C_I(g_j,f_j)_{w}\\
&=&[j,0,j+1,j+2,\cdots,n-1]\cdot(j-1)-[j,j-1,0,j+2,\cdots,n-1]\\
&=&[j-1,0,j,j+1,j+2,\cdots,n-1].
\end{eqnarray*}
Then it is trivial modulo $f_2$ if $j=2$ and modulo $g_{j-1}$ if $j\geq3$. $\Box$

\begin{definition} A \emph{tree} is a connected non-empty graph without circuits.
\end{definition}

A geodesic in a tree is a path without backtracking. The length of
the geodesic from $v$ to $v'$ is called the distance from $v$ to
$v'$, and is denoted by $l(v,v')$.

Fix a vertex $v_0$ of a tree $\Gamma$. For each integer $n\geq0$, let $V_n$ be
the set of vertices $v$ of $\Gamma$ such that $l(v_0,v)=n$.
Then the set of vertices of $\Gamma$ is the union of $V_n$ and $V_i\cap V_j=\emptyset,
 \ i\neq j$. If $v\in V_n$ with $n\geq1$, there is a single vertex $v'\in V_{n-1}$
  from $v_0$ to which $v$ is adjacent.

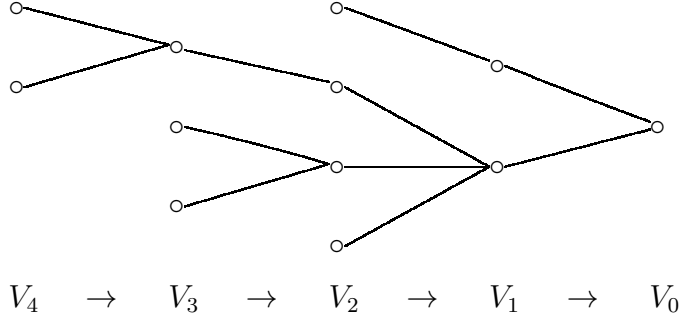
\begin{picture}(100,140)
\put(100,120){$\circ$}\put(100,90){$\circ$}
\put(160,105){$\circ$}\put(160,75){$\circ$}\put(160,45){$\circ$}
\put(220,120){$\circ$}\put(220,90){$\circ$}\put(220,60){$\circ$}\put(220,30){$\circ$}
\put(280,60){$\circ$}\put(280,98){$\circ$}\put(340,75){$\circ$}
\qbezier(340,80)(313,90)(286,101)
\qbezier(280,103)(253,113)(226,123)
\qbezier(340,77)(313,70)(286,63)
\qbezier(280,63)(253,63)(226,63)
\qbezier(280,63)(253,48)(226,33)
\qbezier(280,63)(253,78)(226,93)
\qbezier(166,107)(193,101)(220,95)
\qbezier(166,78)(193,72)(220,64)
\qbezier(166,48)(193,56)(220,64)
\qbezier(106,123)(133,116)(160,109)
\qbezier(106,93)(133,101)(160,109)
\put(100,10){$V_4$}\put(160,10){$V_3$}\put(220,10){$V_2$}\put(280,10){$V_1$}\put(340,10){$V_0$}
\put(129,10){$\rightarrow$}\put(189,10){$\rightarrow$}\put(249,10){$\rightarrow$}\put(309,10){$\rightarrow$}
\end{picture}

We linearly order the set of vertices $V=\bigcup_{n\geq0} V_n$ such
that $v_0$ is the smallest element and for any $v\in V_i, \ v'\in
V_j$, $v<v'$ if $i<j$. Then the partial commutative metabelian Lie
algebra related to the tree $\Gamma$ is defined by:
$$
\mathcal{ML}_{\Gamma}=\mathcal{L}_{(2)}( V  | R),
$$
where
$$
R=\{ [v',v]=0| v'\in V_{n+1}, v\in V_n,  \ v' \ \mbox{and} \ v \ \mbox{are \ adjacent}, n\geq0 \}.
$$

\begin{theorem}\label{t3.7}  The relation set $R$ forms a Gr\"{o}bner-Shirshov basis for the partial
commutative metabelian Lie algebra $ \mathcal{ML}_{\Gamma}$ related
to the tree $\Gamma$.
\end{theorem}
{\bf Proof.} It is obvious that for any $v'\in V_{n+1}$,  there is
only one element $v\in V_n$ such that the relation $[v',v]=0$ lies
in $R$, which means there is no composition in $R$ at all. Thus, $R$
is a Gr\"{o}bner-Shirshov basis automatically. $\Box$

By Theorems \ref{CDML} and \ref{t3.7}, we have the following
corollary.

\begin{corollary} A linear basis of $ \mathcal{ML}_{\Gamma}$ consists of regular
words $v_0v_1\cdots v_n$ $(n\geq0)$ on $V$ satisfying the following
condition: if $v_0>v_i$ $(i\geq1)$, then $l(v_0,v_i)\neq1$.
\end{corollary}

\begin{definition} Let $n$ be a positive integer. An \emph{$n$-cube},
denoted by $Cu_n$, is a graph which the set of vertices
$V_n=\{(\varepsilon_1,\varepsilon_2,\ldots,\varepsilon_n)
\in\mathbb{R}^n|\varepsilon_i=0 \  \mbox{or} \ 1\}$
 and two vertices $\varepsilon=(\varepsilon_1,\varepsilon_2,\ldots,
 \varepsilon_n)$, $\delta=(\delta_1,\delta_2,\ldots,\delta_n)$ are
 adjacent if $\exists  \ i$, such that $\varepsilon_i=\delta_i+1 \ mod \ 2$
 and $\varepsilon_j=\delta_j$ for any $j\neq i$.
\end{definition}

For example, 3-cube and 4-cube are the followings:

\begin{picture}(100,190)

\put(5,150){$Cu_3$:}
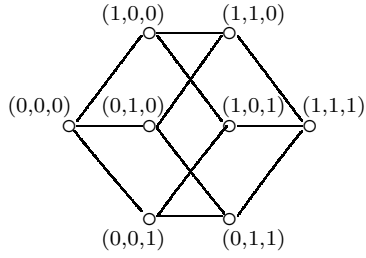
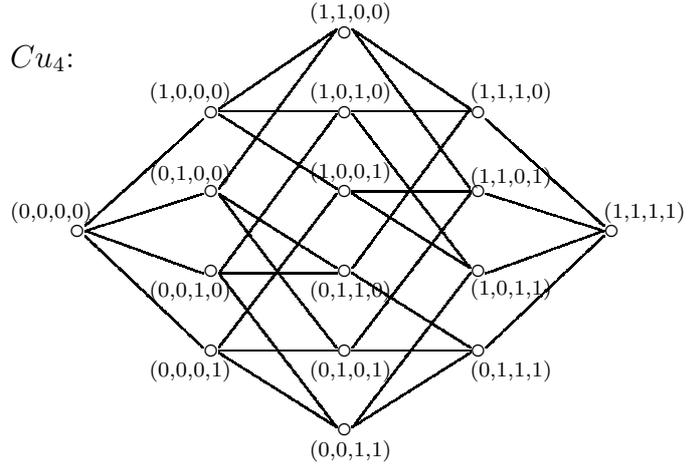
\put(30,85){$\circ$}\put(60,120){$\circ$}
\put(60,85){$\circ$}\put(60,50){$\circ$}
\put(90,120){$\circ$}\put(90,85){$\circ$}
\put(90,50){$\circ$}\put(120,85){$\circ$}
\put(10,95){$_{(0,0,0)}$}\put(45,130){$_{(1,0,0)}$}
\put(45,95){$_{(0,1,0)}$}\put(45,45){$_{(0,0,1)}$}
\put(90,130){$_{(1,1,0)}$}\put(90,95){$_{(1,0,1)}$}
 \put(90,45){$_{(0,1,1)}$}\put(120,95){$_{(1,1,1)}$}

\qbezier(36,90)(48,106)(60,122)
\qbezier(36,88)(48,88)(60,88)
\qbezier(35,86)(47,71)(60,56)
\qbezier(66,123)(78,123)(90,123)
\qbezier(66,122)(78,106)(90,90)
\qbezier(66,87)(78,71)(91,55)
\qbezier(66,54)(78,54)(90,54)
\qbezier(66,54)(79,72)(92,88)
\qbezier(66,88)(78,105)(90,122)
\qbezier(96,122)(108,106)(120,90)
\qbezier(96,88)(108,88)(120,88)
\qbezier(96,54)(108,71)(120,86)

\put(188,150){$Cu_4$:}
\put(210,85){$\circ$}
\put(260,130){$\circ$}\put(260,100){$\circ$}\put(260,70){$\circ$}\put(260,40){$\circ$}
\put(310,160){$\circ$}\put(310,130){$\circ$}\put(310,100){$\circ$}
\put(310,70){$\circ$} \put(310,40){$\circ$} \put(310,10){$\circ$}
\put(360,130){$\circ$}\put(360,100){$\circ$}\put(360,70){$\circ$}\put(360,40){$\circ$}
\put(410,85){$\circ$}

\put(188,95){$_{(0,0,0,0)}$}
\put(240,140){$_{(1,0,0,0)}$}\put(240,109){$_{(0,1,0,0)}$}\put(240,65){$_{(0,0,1,0)}$}\put(240,35){$_{(0,0,0,1)}$}
\put(300,170){$_{(1,1,0,0)}$}\put(300,140){$_{(1,0,1,0)}$}\put(300,110){$_{(1,0,0,1)}$}
\put(300,65){$_{(0,1,1,0)}$}\put(300,35){$_{(0,1,0,1)}$}\put(300,5){$_{(0,0,1,1)}$}
\put(360,140){$_{(1,1,1,0)}$}\put(360,108){$_{(1,1,0,1)}$}\put(360,66){$_{(1,0,1,1)}$}\put(360,35){$_{(0,1,1,1)}$}
\put(410,95){$_{(1,1,1,1)}$}

\qbezier(216,90)(238,110)(260,130)
\qbezier(216,88)(238,95)(260,102)
\qbezier(216,88)(238,80)(260,72)
\qbezier(216,86)(238,65)(260,44)

\qbezier(410,90)(388,110)(366,130)
\qbezier(410,88)(388,95)(366,102)
\qbezier(410,88)(388,80)(366,72)
\qbezier(410,86)(388,65)(366,44)

\qbezier(266,134)(288,149)(310,164)
\qbezier(266,133)(288,133)(310,133)
\qbezier(266,132)(288,118)(310,104)

\qbezier(266,104)(288,134)(310,164)
\qbezier(266,102)(288,88)(310,74)
\qbezier(266,102)(288,72)(310,44)

\qbezier(266,72)(288,102)(310,132)
\qbezier(266,72)(288,72)(310,72)
\qbezier(266,72)(288,43)(310,14)

\qbezier(266,44)(288,73)(310,102)
\qbezier(266,43)(288,43)(310,43)
\qbezier(266,42)(288,28)(310,14)

\qbezier(316,164)(338,149)(360,134)
\qbezier(316,133)(338,133)(360,133)
\qbezier(316,74)(338,104)(360,134)

\qbezier(316,164)(338,134)(360,104)
\qbezier(316,103)(338,103)(360,103)
\qbezier(316,44)(338,73)(360,102)

\qbezier(316,132)(338,103)(360,74)
\qbezier(316,102)(338,88)(360,74)
\qbezier(316,14)(338,43)(360,72)

\qbezier(316,72)(338,58)(360,44)
\qbezier(316,43)(338,43)(360,43)
\qbezier(316,14)(338,28)(360,42)

\end{picture}

We order all vertices lexicographically. The distance of
$\varepsilon$ and $\delta$ is
$d(\varepsilon,\delta)=\sum_{i=1}^n|\varepsilon_i-\delta_i|$. Then
the partial commutative metabelian Lie algebra related to the
$n$-cube $Cu_n$ is defined by:
$$
\mathcal{ML}_{\Gamma}=\mathcal{L}_{(2)}( V_n |\varepsilon\delta=0, \
 d(\varepsilon,\delta)=1).
$$

\begin{theorem}  A Gr\"{o}bner-Shirshov basis $S$ for the partial commutative
metabelian Lie algebra related to 3-cube
$$
 \mathcal{ML}_{Cu_3}=\mathcal{L}_{(2)}(V_3
|\varepsilon\delta, \ d(\varepsilon,\delta)=1, \varepsilon>\delta)
$$
is the union of the following:
\begin{eqnarray*}
&&R_2=\{\lfloor\varepsilon\delta\rfloor \  | \  d(\varepsilon,\delta)=1\},\\
&&R_3=\{\lfloor\varepsilon\delta\rfloor\mu \  |  \ d(\varepsilon,\delta)=2, \ \mu\varepsilon,\mu\delta\in R_1\},\\
&&R_4=\{\lfloor\varepsilon\delta\rfloor\mu\gamma \  |  \ d(\varepsilon,\delta)=3, \ \mu\varepsilon\in R_2,\mu\delta\gamma\in R_3\},\\
&&R_5=\{\lfloor\delta_1\delta_2\rfloor\gamma\langle\mu_1\mu_2\rangle \  |  \ d(\delta_1,\delta_2)=2, \ \gamma\delta_i\mu_i\in R_3, i=1,2\},\\
&&R_5'=\{\lfloor\delta_1\delta_2\rfloor\gamma\mu\mu' \  |  \
d(\delta_1,\delta_2)=2, \ \gamma\delta_1\in R_2, \gamma_2\mu\mu'\in
R_4, d(\mu,\delta_1)\neq1\},
\end{eqnarray*}
where
$\lfloor\varepsilon\delta\rfloor=max\{\varepsilon,\delta\}min\{\varepsilon,\delta\}$.
\end{theorem}

By Algorithm \ref{a1}, we have that a reduced Gr\"{o}bner-Shirshov
basis (it means there is no composition of type I, II, III)  for the
partial commutative metabelian Lie algebra related to $4$-cube $
\mathcal{ML}_{Cu_4}$ consists of 268 relations.

\ \

\noindent{\bf Acknowledgement}: The authors would like to thank
Professor L.A. Bokut for his guidance, useful discussions and
enthusiastic encouragement in writing up this paper.


\begin{thebibliography}{1}


\bibitem{b63}L.A. Bokut, A basis of free polynilpotent Lie algebras, \emph{Algebra Logika}, \textbf{2(4)}(1963),
13-19.

\bibitem{bc}L.A. Bokut and Yuqun Chen, Gr\"{o}bner-Shirshov bases: some new
results, Advances in Algebra and Combinatorics, World Scientific,
2008, 35-56.



\bibitem{bcs}L.A. Bokut, Yuqun Chen and K.P. Shum, Some new results on
Gr\"{o}bner-Shirshov bases, Proceedings of International Conference
on Algebra, Gadjah Mada University, Indonesia, 7-10 October 2010,
World Scientific, to appear. arxiv.org/abs/1102.0449




\bibitem{dkr1}E. Daniyarova, I. Kazatchkov and V. Remeslennikov, Semidomains
and metabelian product of metabelian Lie algebras, \emph{Journal of Mathematical Sciences},
\textbf{131}(6)(2005), 6015-6022.

\bibitem{dkr2}E. Daniyarova, I. Kazachkov and V. Remeslennikov, Algebraic geometry
over free metabelian Lie algebra I: U-algebras and universal classes, \emph{Journal
of Mathematical Sciences}, \textbf{135}(5)(2006), 3292-3310.

\bibitem{dkr3}E. Daniyarova, I. Kazachkov and V. Remeslennikov, Algebraic geometry
over free metabelian Lie algebra II: Finite field case,
\emph{Journal of Mathematical Sciences}, \textbf{135}(5)(2006), 3311-3326.

\bibitem{df}V. Drensky and S. Findik, Inner and outer automorphisms of free
metabelian nilpotent Lie algebras, \emph{Mathematical Physics and Mathematics}, to appear.

\bibitem{f}S. Findik, Normal and normally outer automorphisms of free metabelian
nilpotent Lie algebras, \emph{Serdica Mathematical Journal},
 to appear.

\bibitem{k}V. Kurlin, The Baker-Campbell-Hausdorff formula in the free
metabelian Lie algebra, \emph{Journal of Lie Theory}, \textbf{17}(3)(2007), 525-538.


\bibitem{ta}V.V. Talapov, Algebraically closed metabelian Lie algebras,
\emph{Algebra i Logika}, {\bf 21}(3)(1982), 357-367.

\end{thebibliography}
\end{document}